\documentclass[graybox]{svmult}

\usepackage{mathptmx}       
\usepackage{helvet}         
\usepackage{courier}        
\usepackage{type1cm}        
\usepackage{makeidx}         
\usepackage{graphicx}        
\usepackage{multicol}        
\usepackage[bottom]{footmisc}

\usepackage{amssymb,amsmath}

\makeindex             

\newcommand{\C}{\mathbb{C}}
\newcommand{\Z}{\mathbb{Z}}

\newcommand{\R}{\mathbb{R}}
\newcommand{\Q}{\mathbb{Q}}
\newcommand{\F}{\mathcal{F}}

\newcommand{\Des}{\mathrm{Des}}

\newcommand{\des}{\mathrm{des}}

\newcommand{\cone}{\mathrm{cone}}
\newcommand{\conv}{\mathrm{conv}}

\newcommand{\cd}[2]{\operatorname{cd}({#1},{#2})}
\newcommand{\arc}[2]{\operatorname{arc}({#1},{#2})}
\newcommand{\open}{\text{int}}
\renewcommand{\phi}{\varphi}
\renewcommand{\emptyset}{\varnothing}

\newcommand\commentout[1]{}

\begin{document}

\title*{Unimodality Problems in Ehrhart Theory}

\author{Benjamin Braun}
\institute{Benjamin Braun \at University of Kentucky, Lexington, KY 40506--0027, \email{benjamin.braun@uky.edu}}

\maketitle

\abstract*{
Ehrhart theory is the study of sequences recording the number of integer points in non-negative integral dilates of rational polytopes.
For a given lattice polytope, this sequence is encoded in a finite vector called the Ehrhart $h^*$-vector.
Ehrhart $h^*$-vectors have connections to many areas of mathematics, including commutative algebra and enumerative combinatorics.
In this survey we discuss what is known about unimodality for Ehrhart $h^*$-vectors and highlight open questions and problems.
}


\section{Introduction}

For a sequence $a_0,a_1,a_2,\ldots, a_n$ of real numbers, we define the sequence to be:
\begin{itemize}
\item \emph{unimodal} if there is a $k$ such that $a_i\leq a_{i+1}$ when $i\leq k-1$ and $a_{i}\geq a_{i+1}$ for $i\geq k$.
\item \emph{log-concave} if $a_j^2\geq a_{j-1}a_{j+1}$ for all $1\leq j <n$.
\item \emph{real-rooted} if the generating polynomial $a_0+a_1x+a_2x^2+\cdots +a_nx^n$ has all real zeros.
(We define constant polynomials to be real-rooted.)
\end{itemize}
It is well-known that real-rootedness implies log-concavity, and that if all the $a_i$ are positive then log-concavity implies unimodality.
Further, these implications are strict as shown by the sequences $(1,4,6,4,1)$, $(1,1,1,1,1)$, and $(1,1,2,1,1)$.

Unimodal sequences occur frequently in mathematics, and a wide range of tools from algebra, analysis, probability, and combinatorics are available for establishing unimodality.
One reason combinatorialists are interested in unimodality results is that their proofs often point to interesting and unexpected properties of associated combinatorial, geometric, and algebraic objects.
There are three important survey articles regarding these properties in combinatorics, due to Stanley \cite{stanleylogconcave}, Brenti \cite{brentisurvey}, and Br\"{a}nd\'{e}n \cite{brandensurvey}.
While the surveys by Stanley and Brenti were written over twenty years ago, they remain relevant and central to the study of unimodal sequences in combinatorics.
The recent survey by Br\"{a}nd\'{e}n is intended to complement the other two works; it covers a range of important newer techniques and developments.

One source of interesting sequences of non-negative integers is Ehrhart theory, the study of enumerating integer points in integer dilates of polytopes whose vertices lie in $\Q^n$.
Our focus in this survey is on sequences arising as Ehrhart $h^*$-vectors for \emph{lattice polytopes}, also known as \emph{integer polytopes}, which are polytopes whose vertices lie in $\Z^n$.
While there are many interesting classes of integer polytopes giving rise to unimodal $h^*$-vectors in Ehrhart theory, fascinating questions and conjectures remain open.

Because of the close connection that exists between lattice polytopes and Cohen-Macaulay semigroup algebras, many of the techniques that have been used in attempts to establish unimodality of $h^*$-vectors are algebraic in nature.
The three surveys mentioned above contain a wealth of additional techniques for establishing unimodality, log-concavity, and real-rootedness.
Our hope is that this survey will contribute to increased interest in unimodality problems in Ehrhart theory, and inspire the application of a broader range of tools to these problems.


\section{Lattice Polytopes and Ehrhart Theory}

In this section we review some well-known background regarding Ehrhart theory for lattice polytopes, including a collection of motivating examples for the study of unimodality.
We will frame this discussion in the more general context of the study of integer points in rational cones, of which Ehrhart theory is a special case.

\subsection{Integer Point Transforms and Ehrhart Series}

Further discussion of the ideas in this section, including proofs of stated results, can be found in \cite{BeckRobinsCCD,BrunsHerzogCMR,Hibi,MillerSturmfels,StanleyGreenBook}.

The study of integer points in rational polyhedra is ubiquitous in mathematics, which is not surprising given that this topic is fundamentally about the tension between arithmetic properties of the lattice $\Z^n$ and geometric convexity in $\R^n$.
For any subset $K\subset \R^n$, we define the \emph{integer point transform} of $K$ to be the formal power series
\[
\sigma_K(z)=\sigma_K(z_1,\ldots,z_n):=\sum_{m\in K \cap \Z^n}z^m  \, ,
\]
where we use the shorthand $z^m:=z_1^{m_1}\cdots z_n^{m_n}$.
Given a pointed rational cone 
\[
C=\{x\in \R^n:Ax\geq 0 \} \, ,
\]
where \emph{pointed} requires that $C$ does not contain any non-trivial linear subspace of $\R^n$ and \emph{rational} requires that $A$ have integral (or rational) entries, it is well-known that $\sigma_C(z)$ is equal to a rational function in the field of fractions $\C((z_1,z_2,\ldots,z_n))$.

We will assume for the rest of this survey that all cones $C$ are pointed and rational.
The rationality of $\sigma_C(z)$ for pointed rational cones is a powerful tool in enumerative combinatorics, and the structure underlying many well-known rational generating function identities.
We will soon see that even relatively simple geometric objects like cubes and simplices can have complicated arithmetic properties.

Our key example of a rational cone is the cone over a rational polytope, defined as follows.
Let $P\subset \R^n$ be a \emph{rational} polytope, that is $P$ is the convex hull $P=\conv(v^1,\ldots,v^k)$ of a collection of vectors $v^1,\ldots,v^k\in \Q^n$.
(We will use upper indices to index vectors, and lower indices to index their components.)
Assume that $P$ has dimension $d$, though often we will be interested in situations where $d=n$.
The \emph{cone over $P$} is 
\[
\cone(P):=\text{span}_{\R_{\geq 0}}\{(1,p):p\in P\}\subset \R\times \R^n \, ,
\]
where we consider the new variable to be indexed at $0$, i.e. $x\in \R\times\R^n=\R^{1+n}$ is written $x=(x_0,x_1,\ldots,x_n)$.
We define the \emph{Ehrhart series} of $P$ to be the generating function
\[
E_P(z):=\sigma_{\cone(P)}(z,1,1,\ldots,1) = 1+\sum_{t\in\Z_{\geq 1}}|tP\cap \Z^n|z^t\, ,
\]
where $tP$ is the notation for the $t$-th dilation of $P$.
That the middle and right-hand terms displayed above are equal follows from the observation that the points in $\cone(P)$ with $0$-th coordinate $t$ form a copy of $tP$.

We will frequently encounter polytopes $P$ that are \emph{lattice polytopes}, also known as \emph{integer polytopes} (we will use both terms freely), which are defined by the condition that $P=\conv(v^1,\ldots,v^k)$ for some $v^1,\ldots,v^k\in \Z^n$.
For lattice polytopes, we may combine two well-known theorems due to Ehrhart \cite{Ehrhart} and Stanley \cite{StanleyDecompositions} and conclude that there exist values $h_0^*,\ldots,h_d^*\in \Z_{\geq 0}$ with $h_0^*=1$ such that
\[
E_P(z)=\frac{\sum_{j=0}^dh_j^*z^j}{(1-z)^{d+1}} \, .
\]
We say the polynomial $h^*_P(z):=\sum_{j=0}^dh_j^*z^j$ is the \emph{$h^*$-polynomial} of $P$ and the vector of coefficients is called the \emph{$h^*$-vector} of $P$, denoted $h^*_P$.
That $E_P(z)$ is of this rational form is equivalent to $|tP\cap \Z^n|$ being a polynomial function of $t$ of degree $d$; the non-negativity of the $h^*$-vector is an even stronger property.
Some authors denote $h^*_i$ by $\delta_i$, and often the $h^*$-polynomial is referred to as the $\delta$-polynomial of $P$ \cite{Hibi}.
We follow the convention of using $h^*$ as introduced by Stanley~\cite{stanleyhstarmonotonic} to emphasize the connections between Ehrhart $h^*$-polynomials and $h$-polynomials for simplicial polytopes.

When studying integer-point transforms of cones, one often considers the following properties.
An $n$-dimensional cone is \emph{simplicial} if it is generated by $n$ linearly independent vectors; an $n$-dimensional polytope is a \emph{simplex} if it is the convex hull of $n+1$ vertices.
A pointed rational simplicial cone $C$ is \emph{unimodular} if the minimal length lattice points on the rays generating $C$ form a parallelpiped of volume $\pm 1$ in the lattice obtained by intersecting the affine span of $C$ with $\Z^n$; when $C$ is an $n$-dimensional cone in $\R^n$, this lattice is $\Z^n$ itself.
A lattice simplex $P$ is \emph{unimodular} if $\cone(P)$ is unimodular.
A \emph{triangulation} of an $n$-dimensional cone $C$ (resp. a $d$-dimensional polytope $P$) is a finite collection $T$ of simplicial $n$-cones (resp. $d$-simplices) such that $C$ (resp. $P$) is the union of the elements of $T$ and for any $\Delta_1,\Delta_2\in T$, $\Delta_1\cap \Delta_2$ is a face common to both $\Delta_1$ and $\Delta_2$.
The triangulation is \emph{unimodular} if every cone (resp. simplex) in the triangulation is unimodular.

Due to their role in triangulations and the following lemma, simplicial cones are fundamental objects in the study of integer point transforms.
\begin{lemma}\label{lemma:simpcone}
Let $C = \{\sum_{ j=0 }^n \lambda_j v^j:\lambda_j\in \R_{\geq 0}\}$ be a simplicial cone in $\R^{ 1+n }$ with linearly independent generators $v^0, v^1, \dots, v^n \in \Z^{ 1+n }$. Then
\[
  \sigma_C (z_0, z_1, \dots, z_n) = \frac{ \sigma_{ \Pi } (z_0, z_1, \dots, z_n) }{ \prod_{ j=0 }^n \left( 1 - z^{ v^j } \right) }
\]
where $\Pi := \{\sum_{ j=0 }^n \lambda_j v^j: \lambda_j\in [0,1)\}$.
Furthermore, the number of integer points in $\Pi$ (and thus the number of monomials in $\sigma_{ \Pi } (z_0, z_1, \dots, z_n)$) is given by the determinant of the matrix with column vectors $v^0, v^1, \dots, v^n$.
\end{lemma}
We refer to $\Pi$ as the \emph{fundamental parallelepiped} of $C$; note that it depends on the choice of generators of~$C$.
When $C=\cone(P)$ for a lattice simplex $P$ and Lemma~\ref{lemma:simpcone} is specialized to the case of the Ehrhart series $E_P(z)$, one obtains the following geometric interpretation of the $h^*$-vector of the simplex $P$.
\begin{lemma}\label{lemma:h*simplex}
For any lattice simplex $P\subseteq \R^n$, with $\Pi$ the fundamental parallelepiped of $\cone(P)$, the $h^*$-vector of $P$ satisfies
\[
h_i(P)=\#\{m\in \Pi\cap \Z^{1+n}:m_0=i\} \, .
\]
\end{lemma}

The $h^*$-vector of a lattice polytope $P$ is a fascinating partial invariant; obtaining a general understanding of $h^*$-vectors of lattice polytopes and their geometric/combinatorial implications is currently of great interest.  
A large amount of recent research \cite{BatyrevDualPolyhedra,BeckHosten,BeyHenkWills,BraunEhrhartFormulaReflexivePolytopes,HibiDualPolytopes,HaaseMelnikov,MustataPayne,Payne} has focused on the class of reflexive lattice polytopes.
Recall that the \emph{dual polytope} of a polytope $P$ containing the origin in its interior is 
\[
P^\Delta := \left\{x\in \mathbb{R}^{n}: x\cdot p \leq 1 \text{ for all } p\in P\right\} \, .
\]
An $n$-dimensional lattice polytope $P\subset \mathbb{R}^{n}$ is \emph{reflexive} if $P^\Delta$ is also a lattice polytope.
For example, crosspolytopes and $\pm 1$-cubes (defined below) form dual pairs of reflexives.
Reflexive polytopes were named by Batyrev \cite{BatyrevDualPolyhedra} in the context of mirror symmetry and theoretical physics, and independently discovered in the context of Ehrhart theory by Hibi \cite{HibiDualPolytopes}, who proved the following theorem.

\begin{theorem}[Hibi, \cite{HibiDualPolytopes}] 
A lattice polytope $P\subset \R^n$ is reflexive (up to unimodular equivalence) if and only if $P$ is $n$-dimensional, $P$ contains an interior lattice point, and $h^*_P$ satisfies $h_i^*=h_{n-i}^*$.
\end{theorem}

It is known from work of Lagarias and Ziegler \cite{lagariasziegler} that there are only finitely many reflexive polytopes (up to unimodular equivalence) in each dimension, with one reflexive in dimension one, $16$ in dimension two, $4,319$ in dimension three, and $473,800,776$ in dimension four according to computations by Kreuzer and Skarke~\cite{KreuzerSkarke00}.
The number of five-and-higher-dimensional reflexives is unknown.
In addition to Ehrhart theory, reflexive polytopes play important roles in combinatorial mirror symmetry, toric geometry, and the theory of error-correcting codes and random walks.

Reflexive polytopes are a special case of the more general family of Gorenstein polytopes.
A lattice polytope $P\subset \R^n$ is \emph{Gorenstein of index $k$} if $kP$ is a reflexive polytope, up to translation by an element of $\Z^n$.
More generally, a pointed rational cone $C$ is \emph{Gorenstein} if there exists an integer point $c$ such that
\[
c+\left( \Z^{1+n}\cap C \right) = \Z^{1+n}\cap \open(C) 
\]
where $\open(C)$ denotes the interior of the cone $C$.
This definition of a Gorenstein cone is motivated by results of Stanley and Danilov for Gorenstein semigroup algebras, see \cite[Section 6.3]{BrunsHerzogCMR} for a textbook exposition.
It can be shown that $P$ is Gorenstein if and only if $\cone(P)$ is Gorenstein.
Gorenstein polytopes may also be identified using the following extension of Hibi's theorem, which is a consequence of the Gorenstein property for semigroup algebras.
\begin{theorem}[Stanley, \cite{stanleyhilbertgraded}]
A $d$-dimensional lattice polytope $P$ with 
\[
h^*_P=(h_1^*,\ldots,h_s^*,0,0,\ldots,0)\in \Z^d \, ,
\]
where $h^*_s\neq 0$, is Gorenstein if and only if $h^*_i=h^*_{s-i}$ for all $i$.
\end{theorem}


\subsection{Examples}

We now consider a range of examples of $h^*$-vectors along with polytopes that realize them.
Throughout this subsection we will denote the $i$-th standard basis vector in $\R^n$ by $e_i$.
We will begin with examples of polytopes that have well-behaved $h^*$-vectors.

\begin{example}
The $h^*$-polynomial $h^*_P(z)=1$ is unique to simplices that are unimodularly equivalent to $\Delta_n:=\conv(0,e_1,\ldots,e_n)$, the \emph{standard unimodular simplex}.  
One can easily show that $\displaystyle E_{\Delta_n}(z)=\frac{1}{(1-z)^{n+1}}$ using Lemma~\ref{lemma:simpcone}.
\end{example}

\begin{example}\label{ex:reflexivesimplex}
The polynomial $h^*_P(z)=1+z+z^2+\cdots +z^d$ is the $h^*$-polynomial of the reflexive simplex of minimal volume in $\R^d$, defined by
\[
S_d:=\conv \{e_1,e_2,\ldots,e_d,-\sum_{i=1}^de_i  \} \, .
\]
This computation is again a straightforward exercise using Lemma~\ref{lemma:simpcone}.
Note that the finite geometric series is unimodal but is not real-rooted.
\end{example}

\begin{example}
The polynomial $(1+z)^d$ is the $h^*$-polynomial for many different lattice polytopes.
Two examples are the \emph{crosspolytopes}
\[
C_d:=\conv\{e_1,\ldots,e_d,-e_1,\ldots,-e_d\}
\]
and the family of reflexive simplices given by 
\[
T_d:=\conv\{-\sum_{i=1}^de_i, e_1-\sum_{i=2}^de_i, e_2-\sum_{i=3}^de_i, \cdots, e_{d-1}-e_d, e_d\} \, .
\]
The $h^*$-polynomial for $C_d$ can be computed in many different ways, see \cite[Chapter 2]{BeckRobinsCCD} for one example.
The derivation of $h_{T_d}^*(z)$ can be done using either Lemma~\ref{lemma:simpcone} or by recognizing $T_d$ as an example of an affine free sum of intervals of length $2$ and applying techniques from \cite{beckjayawantmcallister}.

It is immediate that the polynomial $(1+z)^d$ is real-rooted, hence the binomial coefficients are both log-concave and unimodal.
Many additional proofs of unimodality for this polynomial are known.
Observe that the symmetry of the binomial coefficients implies that both $C_d$ and $T_d$ are reflexive; the dual polytope to $C_d$ is the cube $[-1,1]^d$.
\end{example}

\begin{example}
Let $S_n$ denote the symmetric group on $n$ objects.
For an element $\pi\in S_n$, the \emph{descent set} of $\pi$ is $\Des(\pi):=\{i\in [n-1]:\pi(i)>\pi(i+1)\}$ and the \emph{descent statistic} for $\pi$ is $\des(\pi):=|\Des(\pi)|$.
The \emph{$n$-th Eulerian polynomial} is $A_n(z):=\sum_{\pi\in S_n}z^{\des(\pi)}$.
The coefficients of the Eulerian polynomials are called the Eulerian numbers.
Eulerian polynomials are known to have symmetric coefficients and to be real-rooted.
This family of polynomials has inspired many subsequent developments regarding real-rootedness, log-concavity, and unimodality in enumerative combinatorics.

The \emph{unit cube} $[0,1]^n$ satisfies $h^*_{[0,1]^n}(z)=A_{n-1}(z)$.
While there are many proofs of this, a geometric proof involving triangulations can be found in \cite{BeckBraunEulerMahonian}.
Another family of polytopes with Eulerian polynomials as $h^*$-polynomials are the simplices
\[
L_d:=\conv\left\{
\begin{array}{l}
e_d, (d-1)e_{d-1}+de_d, \\
(d-2)e_{d-2}+ (d-1)e_{d-1}+de_d,\ldots, \\
2e_2+\cdots +de_d, e_1+2e_2+\cdots +de_d 
\end{array}\right\} \, .
\]
That $L_d$ is a $(d-1)$-dimensional simplex in $\R^d$ satisfying $h^*_{L_d}(z)=A_{d-1}(z)$ is implicit in work of Corteel, Lee, and Savage \cite{CLS2005}.
Observe that the symmetry of the Eulerian numbers implies that $[0,1]^n$ is Gorenstein of index 2 (a straightforward verification) and $L_d$ is reflexive.
\end{example}

\begin{example}
There is an analogous definition of the descent statistic for the hyperoctahedral group, i.e. the group $B_n$ of signed permutations, which gives rise to the \emph{type B Eulerian polynomial} $\sum_{\pi\in B_n}z^{\des(\pi)}$.
The type B Eulerian polynomials are also known to be real-rooted with symmetric coefficients.
The cube $[-1,1]^n$ satisfies 
\[
h^*_{[-1,1]^n}(z)=\sum_{\pi\in B_n}z^{\des(\pi)} \, ,
\]
and a geometric proof can be found in \cite{BeckBraunEulerMahonian}.
As with the Eulerian numbers, the type B Eulerian polynomial can also be obtained as the $h^*$-polynomial of a simplex, namely
\[
L_n^2:=\{x\in \R^n|0\leq \frac{x_1}{2}\leq \frac{x_2}{4}\leq \frac{x_3}{6}\leq \cdots \leq \frac{x_n}{2n}\leq 1 \} \, .
\]
This was first shown by Savage and Schuster \cite{savageschuster}.

Again, the symmetry of the type B Eulerian polynomials implies that both the $\pm 1$-cube and $L_n^2$ are reflexive.
\end{example}

\begin{example}
Another polynomial with symmetric log-concave coefficients is the \emph{Narayana polynomial} $N_n(z):=\sum_{j=1}^n N(n,j)z^{j-1}$, where $N(n,j)=\frac{1}{n}\binom{n}{j}\binom{n}{j-1}$.
It is well-known \cite[Chapter 6]{StanleyVol2} that $N_n(1)$ is equal to the $n$-th Catalan number, and that $N(n,j)$ arises in many situations involving natural refinements of Catalan objects.
Let $P_n\subset \R^n$ denote the convex hull of the positive simple roots of type A together with the origin, i.e.
\[
P_n :=\conv\{0,e_i-e_j:n\geq i>j\geq 1\} \, .
\]
It follows from work of Gelfand, Graev, and Postnikov \cite{GelfandGraevPostnikov} combined with techniques due to Stanley as described in \cite[Solution to Exercise 6.31(b)]{StanleyVol2} that $h^*_{P_n}(z)=N_n(z)$, hence $P_n$ is Gorenstein.
\end{example}

In case the preceding examples have painted an overly-optimistic portrait of the world of Ehrhart $h^*$-vectors, we end this subsection with two examples of situations where $h^*$-polynomials behave in pathological ways.

\begin{example}
Let $m,d,k\in \Z_{\geq 0}$ satisfy $m\geq 1$, $d\geq 2$, and $1\leq k\leq \lfloor (d+2)/2\rfloor$.
For each such triple of integers, Higashitani \cite{higashitanicounterexamples} has constructed a $d$-dimensional lattice simplex with $h^*$-polynomial $1+mz^k$.
The construction is relatively simple; verifying the $h^*$-polynomial requires only the use of Lemma~\ref{lemma:simpcone}.
\end{example}

\begin{example}
Payne~\cite{Payne} gives an explicit construction showing that for any positive integers $m$ and $n$, there exists a reflexive polytope $P$ and indices 
\[
i_1<j_1<i_2<j_2<\cdots <i_m<j_m<i_{m+1}
\]
such that the $h^*$-vector for $P$ satisfies
\[
h^*_{i_\ell}-h^*_{j_\ell}\geq n \, \text{ and } \,  h^*_{i_{\ell+1}}-h^*_{j_\ell}\geq n \, .
\]
Thus, reflexive polytopes exist whose $h^*$-vectors violate unimodality arbitrarily badly with arbitrary length.
For all dimensions greater than or equal to $6$, Payne also constructs explicit examples of reflexive simplices with non-unimodal $h^*$-vectors.
\end{example}


\subsection{Integral Closure, Normality, and the Integer Decomposition Property (IDP)}

Many of the connections between polytopes and unimodality involve the notions of integral closure and normality.
See Bruns and Gubeladze \cite{brunsgubeladzebook} for a complete discussion of the distinction between these two ideas.
A lattice polytope $P$ in $\R^n$ is \emph{normal} if
\[
\text{span}_{\Z_{\geq 0}}\{(1,P)\cap \Z^{1+n}\}=\cone(P)\cap \text{span}_{\Z}\{(1,P)\cap \Z^{1+n}\} \, .
\]
$P$ is said to be \emph{integrally closed}, or to satisfy the \emph{integer decomposition property (IDP)}, if 
\[
\text{span}_{\Z_{\geq 0}}\{(1,P)\cap \Z^{1+n}\}=\cone(P) \cap \Z^{1+n}\, .
\]
Thus, normality for lattice polytopes involves integer points in $\cone(P)$ that are in the $\Z$-linear span of integer points in $(1,P)$, while integral closure/IDP for lattice polytopes involves all integer points in $\cone(P)$.
It follows directly from the definitions that IDP implies normal.

In the research literature on this topic, it is common for authors to assume that $\Z^{1+n}=\text{span}_{\Z}\{(1,P)\cap \Z^{1+n}\}$, i.e. to consider the right-hand lattice as the lattice of reference for $P$, and it is also common for authors to use these terms as synonyms.
As a result, the terms ``normal'' and ``integrally closed'' are potential causes of confusion in the geometric combinatorics literature, see~\cite[Remark 0.1]{coxhaasehibihigashitani} for examples where this arises.
For polytopes where the two possible lattices of reference are not equal, one must be careful to specify which is being used for lattice point enumeration; we will assume throughout that $\Z^{1+n}$ is being considered unless specified otherwise.

Several authors~\cite{coxhaasehibihigashitani} have proposed to replace the name ``integrally closed'' with ``IDP'' for two reasons: first, to alleviate the existing confusion in the literature between integrally closed and normal polytopes, and second to align with the standard use of IDP as the nomenclature in integer programming and optimization.
For these reasons, in this survey we will use the term IDP when discussing integrally closed polytopes.


\subsection{Connections to Commutative Algebra}

The definitions of integral closure and normality are closely related to connections between Ehrhart theory and Hilbert series for finitely generated graded algebras.
Given a finitely generated $\C$-algebra $A=\oplus_{i=0}^\infty A_i$, we say $A$ is \emph{standard} if $A$ is generated by $A_1$ and \emph{semistandard} if $A$ is integral over the subalgebra generated by $A_1$.
Given a lattice polytope $P$ in $\R^n$, the \emph{semigroup algebra associated to $P$} is 
\[
\C[P]:=\C[x^v:v\in \cone(P)\cap \Z^{1+n}] \, .
\]
We grade $C[P]$ by $\deg(x_0^{m_0}\cdots x_n^{m_n})=m_0$, i.e. by the ``height'' of the exponent vector in $\cone(P)$.
With this grading, $\C[P]$ is a semistandard semigroup algebra, where the algebra generated by $\C[P]_1$ is $\C[x^v:v\in (1,P)\cap \Z^{1+n}]$, i.e. the semigroup algebra generated by integer points in $(1,P)$.
The role played by the IDP property is captured by the following proposition, whose proof is immediate.

\begin{proposition}
A lattice polytope $P$ is IDP if and only if $\C[P]$ is standard.
\end{proposition}

\begin{example}
For the standard unimodular simplex $\Delta_n$, we have
\[
\C[\Delta_n]=\C[x_0,x_0x_1,x_0x_2,\ldots,x_0x_n] \, ,
\]
where every generator of the algebra is of degree $1$.
Thus, $\C[\Delta_n]$ is isomorphic to a polynomial ring in $n+1$ variables.
\end{example}

\begin{example}
The interval $P=[0,2]$ is IDP, thus we have
\[
\C[P]=\C[x_0,x_0x_1,x_0x_1^2] \, .
\]
Note that this algebra is not a polynomial ring, due to the existence of relations among the generators, e.g. $(x_0)(x_0x_1^2)-(x_0x_1)^2=0$.
However, $\C[P]$ is integral over the subalgebra $\C[x_0,x_0x_1^2]$, which is isomorphic to a polynomial ring.
\end{example}

\begin{remark}
While $\C[P]$ is always integral over $\C[x^v:v\in (1,P)\cap \Z^{1+n}]$, the algebra generated by $\C[P]_1$, it is not necessarily the normalization of the latter algebra.
In general, the normalization of $\C[x^v:v\in (1,P)\cap \Z^{1+n}]$ is $\C[x^v:v\in \cone(P)\cap \text{span}_{\Z}\{(1,P)\cap \Z^{1+n}]$.
If $\Z^{1+n}=\text{span}_{\Z}\{(1,P)\cap \Z^{1+n}\}$, then these two algebras are the same.
\end{remark}
Recall that the Hilbert series of a graded algebra $A$ is the generating function
\[
H(A;z):=\sum_{i=0}^\infty \dim_{\C}(A_i) z^i \, .
\]
It is a straightforward consequence of the definition of $\C[P]$ that
\[
H(\C[P];z)=E_P(z) \, .
\]

By a theorem of Hochster \cite{Hochster}, $\C[P]$ is a Cohen-Macaulay integral domain.
Further, $\C[P]$ always admits a linear system of parameters \cite[Lemma 33.7]{Hibi}.
Thus, if $\dim(P)=d$, there exist algebraically independent elements $\theta_1,\ldots,\theta_{d+1}\in\C[P]_1$ such that $\C[P]$ is a finitely generated free module over the subalgebra $\C[\theta_1,\ldots,\theta_{d+1}]$.
Many important properties of Ehrhart series were first established using this property of $\C[P]$, e.g. the non-negativity of $h^*_P(z)$ \cite{StanleyGreenBook}.
In this context, the symmetry of the $h^*$-vectors for reflexive polytopes implies that $\C[P]$ are Gorenstein when $P$ is reflexive.
Gorenstein algebras form an important subclass of Cohen-Macaulay algebras, and the connection with our previous definition of Gorenstein lattice polytope is that $P$ is Gorenstein if and only if $\C[P]$ is Gorenstein.

One of the open problems of interest in Ehrhart theory is a special case of the following conjecture due to Stanley \cite[Conjecture 4(a)]{stanleylogconcave}, also appearing in Hibi \cite[Conjecture 1.5]{hibiflawless}.
\begin{conjecture}[Stanley, \cite{stanleylogconcave}]\label{conj:stanleygor}
For a standard graded Gorenstein integral domain $A$ of Krull dimension $d$, with
\[
H(A;z)=\frac{h_0+h_1z+\cdots +h_sz^s}{(1-z)^d} \, ,
\]
the sequence $h_0,h_1,\ldots,h_s$ is unimodal.
\end{conjecture}
Stanley's original conjecture was that standard graded Cohen-Macaulay integral domains have unimodal $h$-vectors\footnote{The original, published version of this article incorrectly stated that counterexamples are known to this conjecture. However, this counterexample is still open as of October 2017. See the appendix at the end of this paper for clarification.}.


\section{Unimodality for Polytopes with the Integer Decomposition Property}\label{sec:intclosed}

When Conjecture~\ref{conj:stanleygor} is restricted to algebras of the form $\C[P]$, we encounter our first major open problem regarding unimodality for Ehrhart $h^*$-vectors.
This problem was proposed by Hibi and Ohsugi in the context of normal polytopes~\cite{hibiohsugiconj}.
\begin{conjecture}[Hibi and Ohsugi, \cite{hibiohsugiconj}]\label{conj:hibiohsugi}
If $P$ is Gorenstein and IDP, then $h^*_P$ is unimodal.
\end{conjecture}
Our second major open problem is a relaxation of Conjecture~\ref{conj:hibiohsugi} due to Schepers and Van Langenhoven~\cite{schepersvanl}, arising as a subconjecture of Stanley's original unimodality conjecture for standard graded Cohen-Macaulay integral domains.
\begin{question}[Schepers and Van Langenhoven, \cite{schepersvanl}]\label{question:svl}
If $P$ is IDP, is it true that $h^*_P$ is unimodal?
\end{question}
These are the outstanding open problems regarding unimodality in Ehrhart theory.
A first indication that these problems might be difficult is that they both fail when considering the real-rooted property instead of unimodality, as demonstrated by the simplex $S_d$ in Example~\ref{ex:reflexivesimplex}.

There are two other well-known results regarding $h^*$-polynomials that further motivate this line of investigation.
First, Conjecture~\ref{conj:hibiohsugi} is rooted in an earlier conjecture due to Hibi that all reflexive polytopes have unimodal $h^*$-vectors.
This was shown to be false in even dimensions by Musta{\c{t}}{\v{a}} and Payne~\cite{MustataPayne} and in all dimensions by Payne~\cite{Payne}.
While their counterexamples arise as reflexive simplices, these simplices are not normal, hence not IDP.
Second, Stanley had conjectured that the $h^*$-vector of the Birkhoff polytope, i.e. the polytope of doubly-stochastic matrices, is unimodal -- this was proved to be true by Athanasiadis~\cite{athanasiadisbirkhoff}.
The Birkhoff polytope is both Gorenstein and IDP, and remains a source of interesting problems regarding unimodality~\cite{davissymmetricbirkhoff}.
A key property of the Birkhoff polytope that Athanasiadis used is that it admits a regular unimodular triangulation.
Subsequent work in this area exploited this property further, so we will elaborate on such triangulations next.

\subsection{Regular Unimodular Triangulations}

We will assume that all triangulations of a lattice polytope $P\subset \R^n$ have a vertex set contained in $\Z^n$.
If $P$ admits a unimodular triangulation (or even a covering by unimodular subsimplices), then $P$ is IDP because $\cone(P)$ is a union of unimodular cones with lattice-point generators of degree $1$.
A triangulation of $P$ is regular if it arises as the projection of the lower hull of a lifting of the lattice points of $P$ into $\R^{1+n}$.
The following theorem due to Athanasiadis, which had been previously discovered by Hibi and Stanley but never published, places strong restrictions on the ``tail'' of the $h^*$-vector for any polytope admitting a regular unimodular triangulation.
\begin{theorem}[Athanasiadis, \cite{athanasiadisstable}]\label{thm:athanasiadistail}
Let $P$ be a $d$-dimensional lattice polytope with $h^*_P=(h_0^*,\ldots,h_d^*)$.
If $P$ admits a regular unimodular triangulation, then $h_i^*\geq h^*_{d-i+1}$ for $1\leq i\leq \lfloor(d+1)/2\rfloor$, 
\[
h^*_{\lfloor(d+1)/2\rfloor}\geq \cdots \geq h^*_{d-1}\geq h^*_d
\]
and
\[
h_i^*\leq \binom{h^*_1+i-1}{i}
\]
for $0\leq i\leq d$.
\end{theorem}
As a corollary of Theorem~\ref{thm:athanasiadistail}, it follows that if $P$ is reflexive and admits a regular unimodular triangulation, then $h^*_P$ is unimodal.
Athanasiadis's proof of Theorem~\ref{thm:athanasiadistail} relies on lemmas due to Stanley and Kalai on the Hilbert $h$-vectors of Cohen-Macaulay subcomplexes of the boundary complex of a simplicial polytope.

Athanasiadis's result was further extended by Bruns and R\"{o}mer \cite{BrunsRomer}, who proved that all $h^*$-vectors of IDP Gorenstein polytopes are $h^*$-vectors for IDP reflexive polytopes.
\begin{theorem}[Bruns and R\"{o}mer, \cite{BrunsRomer}]
If $P$ is Gorenstein and IDP, then $h^*_P$ is the $h^*$-vector of an IDP reflexive polytope.
Further, if $P$ admits a regular unimodular triangulation, then there exists a simplicial polytope $Q$ such that $h^*_P$ is the $h$-vector of $Q$, and hence $h^*_P$ is unimodal as a consequence of the $g$-theorem.
\end{theorem}
Thus, the presence of a regular unimodular triangulation forces strong restrictions on the $h^*$-vector of a Gorenstein polytope.

In general, the $h^*$-vectors of lattice polytopes satisfy an array of interesting linear inequalities.
Stapledon~\cite{stapledoninequalities} found combinatorial proofs for many of these, including analogues of Theorem~\ref{thm:athanasiadistail}.
In order to establish these inequalities, for a $d$-dimensional lattice polytope $P$ with $\deg(h^*_P(z))=s$ Stapledon introduces the polynomial 
\[
\overline{h^*_P}(z):=(1+z+z^2+\cdots +z^{d-s})h^*_P(z)
\]
and proves that $\overline{h^*_P}(z)=a(z)+z^{d+1-s}b(z)$ for symmetric polynomials $a(z)$ and $b(z)$ of degree $d$ and $s-1$, respectively, with integer coefficients.
Many of the inequalities for $h^*_P$ that Stapledon produces are consequences of properties of $\overline{h^*_P}(z)$, $a(z)$, and $b(z)$; the following theorem is a representative example of the role played by these polynomials.
\begin{theorem}[Stapledon, \cite{stapledoninequalities}]
If $P$ is a $d$-dimensional lattice polytope such that the boundary of $P$ admits a regular unimodular triangulation, then $a(t)$ is unimodal.
\end{theorem}
$P$ is reflexive if and only if $h^*_P(z)=a(z)$.
Because a reflexive polytope $P$ admits a regular unimodular triangulation if and only if the boundary of $P$ admits a regular unimodular triangulation, we see that this is an extension of the result of Athanasiadis.

While regular unimodular triangulations are a powerful tool \cite{deloerarambausantos}, not all lattice polytopes admit such decompositions, even when $P$ is normal \cite{hibiohsugi01}.
Thus, we turn our attention to approaches to establishing unimodality that do not rely on this property.


\subsection{Lefschetz Elements}

Let $A=\oplus_{i=0}^sA_i$ be a finitely generated graded algebra of Krull dimension zero.
A useful technique for establishing unimodality for sequences arising as 
\[
\dim_\C(A_0),\dim_\C(A_1),\ldots,\dim_\C(A_s)
\]
is to find a Lefschetz element in $A$, defined as follows.
A linear form $l \in A_1$ is called a {\em weak Lefschetz element} if the multiplication map
\[
\times l : A_i \to A_{i+1}
\]
has maximal rank, that is, is either injective or surjective, for each $i$.
By Remark 3.3 of \cite{HarimaLefschetz}, if $A$ has a weak Lefschetz element, then the Hilbert series $H(A;z)$ has unimodal coefficients. 
For a $d$-dimensional lattice polytope $P$, the general theory of Cohen-Macaulay algebras implies that if $\theta_1,\ldots,\theta_{d+1}\in\C[P]_1$ is a linear system of parameters for $\C[P]$, then 
\[
h^*_P(z)=H(\C[P]/(\theta_1,\ldots,\theta_{d+1});z) \, ,
\]
i.e. $h^*_i=\dim_\C[\C[P]/(\theta_1,\ldots,\theta_{d+1})]_i$.
Thus, the existence of a Lefschetz element in $\C[P]/(\theta_1,\ldots,\theta_{d+1})$ implies $h^*$-unimodality for $P$.

Several important recent results and counterexamples in Ehrhart theory \cite{hibihermite,higashitanicounterexamples} have involved only simplices.
Also, as shown in the examples given previously, many $h^*$-polynomials can be realized using simplices.
This suggests that lattice simplices are richer objects than they first appear.
In the case of a simplex $P$ with vertices $v^0,\ldots,v^n$, we saw in Lemma~\ref{lemma:h*simplex} that $h^*_i$ counts the number of lattice points of degree $i$ in the fundamental parallelepiped for $P$.
Further, the monomials $x_0x^{v^0},\ldots,x_0x^{v^n}$ form a linear system of parameters in $\C[P]$.
This motivates the study of the zero-(Krull)-dimensional algebra
\[
R_P:=\C[P] / (x_0x^{v^0},\ldots,x_0x^{v^n})
\]
graded by the exponent on $x_0$, which has a basis given by the lattice points in the fundamental parallelepiped for $P$.
Because of the explicit connection between $h^*$-vectors, Hilbert series, and $R_P$, the author and Davis \cite{BraunDavisReflexive} have initiated the investigation of $h^*$-unimodality for IDP reflexive simplices using Lefschetz techniques.

While experimental data suggests that a weak Lefschetz element exists for many of the algebras $R_{P}$ when $P$ is an IDP reflexive simplex, such an element need not exist.
\begin{theorem}[Braun and Davis, \cite{BraunDavisReflexive}]\label{theorem:nolefschetz}
Denote by $Lef_d$ the convex hull of the vectors
\[
e_1, \, \ldots, \, e_d, \, -de_1-\sum_{k=2}^d(d+1)e_k \, .
\]
Then 
\[
h^*_{Lef_d}(z)=1+(d+2)z+(d+2)z^2+\cdots + (d+2)z^{d-1}+z^d
\]
and $R_{Lef_d}$ does not have a weak Lefschetz element.
\end{theorem}
An important observation is that for other choices of linear system of parameters, the quotient of $\C[Lef_d]$ by such a system does have a weak Lefschetz element.
However, the failure of the weak Lefschetz approach when using the vertices of $P$ as a system of parameters complicates the situation, as most of the techniques in Ehrhart theory for studying $h^*$-vectors of simplices correspond to studying the algebra $R_P$ with this choice of system of parameters.
Thus, we pose the following questions.
\begin{question}
For which IDP reflexive simplices does $R_P$ admit a Lefschetz element?
\end{question}
\begin{question}
For $P$ an IDP reflexive simplex, is there a canonical choice of linear system of parameters such that $\C[P]/(\theta_1,\ldots,\theta_{d+1})$ admits a Lefschetz element?
\end{question}

Some evidence in support of $h^*$-unimodality for arbitrary IDP reflexive simplices is that this property is preserved under affine free sums.
Suppose $P\subset \R^n$ and $Q \subset \R^m$ are full-dimensional simplices with $0 \in P$ and $\{v^0, \ldots, v^m\}$ denoting the vertices of $Q$. 
Then for each $i = 0, 1, \ldots, m$ we define
\[
P *_i Q := \conv{(P \times 0) \cup (0 \times Q - v^i)} \subseteq \R^{n+m} \, .
\]
The polytope $P*_i Q$ is an example of an \emph{affine free sum} \cite{beckjayawantmcallister}.
\begin{theorem}[Braun and Davis, \cite{BraunDavisReflexive}]
If $P$ and $Q$ are IDP reflexive simplices with $0\in \open(P)$, then so is $P *_i Q$ for each $i$. 
If, in addition, $h^*_P$ and $h^*_Q$ are unimodal, then so is $h^*_{P *_i Q}$.
\end{theorem}

\begin{remark}
We are not aware of examples of IDP reflexive simplices that do not admit a regular unimodular triangulation.
It would be interesting to know if simplices with these properties exist, and if so, to have explicit constructions of infinite families.
\end{remark}


\subsection{Box Unimodal Triangulations}

Another approach to proving $h^*$-unimodality for IDP reflexive polytopes was proposed by Schepers and Van Langenhoven~\cite{schepersvanl}.
For a lattice simplex 
\[
S=\conv(v^0,v^1,\ldots,v^k)\subset\R^n \, ,
\]
define the \emph{open fundamental parallelepiped of $S$} to be
\[
\Pi^\circ:=\left\{\sum_{i=0}^k\lambda_i(1,v^i):\lambda_i\in (0,1)\right\}
\]
and the \emph{box polynomial of $S$} to be
\[
B_S(z):=\sum_{m\in \Pi^\circ\cap \Z^{1+n}}z^{m_0} \, ,
\]
where $B_\emptyset(z):=1$.
\begin{theorem}[Betke and McMullen, \cite{betkemcmullen}]\label{thm:betkemcmullenlink}
Let $P$ be a reflexive polytope and let $T$ be a triangulation of the boundary of $P$.
Then
\[
h^*_P(z)=\sum_{F \, \text{ face of } \, T}B_F(z)h_F(z)
\]
where $h_F(z)$ is the $h$-polynomial of the link of $F$ in $T$.
\end{theorem}
It is well-known that if $P$ is $d$-dimensional and $T$ is regular, then the polynomials $h_F(z)$ are symmetric, unimodal, have non-negative integer coefficients, and have degree $d-1-\dim(F)$.
A triangulation $T$ of the boundary of $P$ is called \emph{box unimodal} if $T$ is regular and $B_F(z)=\sum_{i=0}^{\dim(F)}b_iz^i$ is symmetric and unimodal for all non-empty faces of $T$.
It follows from Theorem~\ref{thm:betkemcmullenlink} that if $P$ admits a box unimodal triangulation, then $h^*_P(z)$ is also symmetric and unimodal.
An affirmative answer to the following question would therefore imply $h^*$-unimodality for IDP reflexive polytopes.
\begin{question}[Schepers and Van Langenhoven, \cite{schepersvanl}]
Does the boundary of every IDP reflexive polytope have a box unimodal triangulation?
\end{question}
One piece of supporting evidence for the existence of box unimodal triangulations is that for every lattice simplex $S$, $B_S(z)$ is symmetric.
This is a straightforward consequence of the existence of an involution on the lattice points in $\Pi^\circ$ that sends 
\[
\sum_{i=0}^k\lambda_i(1,v^i)\in \Pi^\circ\cap \Z^{1+n}
\]
to 
\[
\sum_{i=0}^k(1-\lambda_i)(1,v^i)\in \Pi^\circ\cap \Z^{1+n} \, .
\]
However, there are at least three complications that arise when considering box unimodal triangulations.
First, there exist simplices with non-unimodal box polynomials~\cite[Example 4.3]{schepersvanl}.
Hence, if the box unimodal triangulation approach is to eventually succeed, one must show that all such simplices are avoided by at least one regular triangulation of the boundary of every IDP reflexive polytope.
Second, a result of Haase and Melnikov~\cite{HaaseMelnikov} states that every lattice polytope is a face of some reflexive polytope.
Thus, it isn't possible to simply ignore examples of simplices with non-unimodal box polynomials -- one must make arguments for why such simplices do not appear in triangulations of IDP reflexives.
Third, it isn't clear what geometric or arithmetic conditions on lattice simplices lead to unimodality for their box polynomials, which leads to the following questions.
\begin{question}
Which lattice simplices have unimodal box polynomials?
\end{question}
\begin{question}
What are necessary conditions for a sequence to be the coefficients of a box polynomial for a lattice simplex?  Sufficient conditions?
\end{question}

\subsection{Zonotopes}

One of the reasons that Gorenstein lattice polytopes are a focus for questions about unimodality is that there are many techniques for establishing unimodality that only apply to symmetric sequences, e.g. the representation theory of $sl_2$ \cite{stanleylogconcave}.
By omitting the requirement that $P$ be Gorenstein, with Question~\ref{question:svl} Schepers and Van Langenhoven remove both a strong geometric constraint on the lattice polytopes under consideration and many established techniques for demonstrating unimodality.

Partial answers to Question~\ref{question:svl} have been obtained.
We call $P$ a \emph{closed lattice parallelepiped} if there exist linearly independent vectors $v^1,\ldots,v^r\in \Z^n$ such that
\[
P=\left\{\sum_{i=1}^r\lambda_iv^i:\lambda_i\in[0,1]\right\} \, ,
\]
i.e. if $P$ is a Minkowski sum of linearly independent lattice segments.
Lattice parallelepipeds are known to be IDP.
\begin{theorem}[Schepers and Van Langenhoven, \cite{schepersvanl}]
Closed lattice parallelepipeds have unimodal $h^*$-vectors.
\end{theorem}
Lattice parallelepipeds are a subclass of \emph{lattice zonotopes}, which by definition are lattice polytopes arising as Minkowski sums of arbitrary lattice segments \cite{mcmullenzonotopes}.
Lattice zonotopes have a covering by lattice parallelepipeds, hence are IDP.
A natural question is if the $h^*$-vectors of lattice zonotopes are unimodal, which the following theorem answers positively.
\begin{theorem}[Beck, Jochemko, and McCullough, \cite{beckjochemkomccullough,jochemkothesis}]
If $P$ is a lattice zonotope, then $P$ has a unimodal $h^*$-vector.
\end{theorem}

In general, IDP lattice polytopes can behave in counterintuitive ways, as shown by the existence of normal lattice polytopes (which are IDP with respect to the lattice $\text{span}_{\Z}\{(1,P)\cap \Z^{1+n}\}$) that do not admit a regular unimodular triangulation~\cite{hibiohsugi01}.
While it seems reasonable to replace the general problem of $h^*$-unimodality for IDP $P$ by the problem of demonstrating $h^*$-unimodality for special subfamilies of IDP polytopes, even this can create significant challenges as shown in the next subsection.


\subsection{Matroid Polytopes}

Recall that a \emph{matroid} $M$ is a system $\F$ of subsets of $\{1,2,\ldots,n\}$ called \emph{independent sets} such that the following hold:
\begin{itemize}
\item $\emptyset\in\F$
\item if $X\in \F$ and $Y\subseteq X$, then $Y\in\F$
\item if $X,Y\in\F$ and $|X|=|Y|+1$, there exists $x\in X\setminus Y$ such that $Y\cup x\in\F$
\end{itemize}
The \emph{bases} of $M$ are the inclusion-maximal independent sets in $M$, and the set of bases for $M$ is denoted $\mathcal{B}$.
The \emph{matroid polytope} of $M$ is
\[
P(M):=\conv\{e_B:B\in\mathcal{B}\} \, ,
\]
where $e_B:=\sum_{i\in B}e_i$ with $e_i$ denoting the $i$-th standard basis vector in $\R^n$.
The following conjecture due to De Loera, Haws, and K\"{o}ppe is based on extensive computational evidence.
\begin{conjecture}[De Loera, Haws, and K\"{o}ppe, \cite{deloerahawskoeppe}]\label{conj:deloerahawskoppe}
For any matroid $M$, the $h^*$-vector of $P(M)$ is unimodal.
\end{conjecture}
Matroid polytopes have many nice properties, for example Ardila, Benedetti, and Doker proved that they are generalized permutohedra~\cite{ardilabenedettidoker}, a class of polytopes defined by Postnikov~\cite{PostnikovAssPermBeyond}.
A geometric variant of an algebraic conjecture due to White~\cite{whitetoric} regarding the toric ideal of $P(M)$ is the following conjecture of Haws.
\begin{conjecture}[Haws, \cite{hawsthesis}]
For any matroid $M$, the polytope $P(M)$ has a regular unimodular triangulation.
\end{conjecture}
Haws also proposed the following weaker conjecture.
\begin{conjecture}[Haws, \cite{hawsthesis}]
For any matroid $M$, the polytope $P(M)$ is the union of its unimodular subsimplices.
\end{conjecture}

We will focus in this subsection on the special case of Conjecture~\ref{conj:deloerahawskoppe} when $M$ is a uniform matroid.
Recall that the \emph{uniform matroid of rank $k$ on $n$ elements}, denoted $U^{k,n}$, is the matroid with bases given by the set of all $k$-element subsets of $\{1,2,\ldots,n\}$.
The polytope $P(U^{k,n})$ is known in the geometric combinatorics literature as the \emph{$(n,k)$-hypersimplex}, denoted $\Delta_{n,k}$, and we will use this notation throughout.
Several results regarding lattice point enumeration for hypersimplices are known~\cite{lampostnikovalcoved,lihypersimplices}.
The unimodality of $h^*_{\Delta_{n,2}}$ follows from work of Katzman~\cite{katzman}, but for all $k\geq 3$ the unimodality of the $h^*$-vectors for hypersimplices remains an open problem.
While not all hypersimplices are Gorenstein, $\Delta_{n,k}$ always admits a regular unimodular triangulation, which we will refer to as the \emph{circuit triangulation}.
While different descriptions of this triangulation have been found, Lam and Postnikov proved that these descriptions all correspond to the same triangulation~\cite{lampostnikovalcoved}.

The author and Solus~\cite{BraunSolusHyper} have introduced the following polytopes as part of an attempt to prove that $\Delta_{n,k}$ is $h^*$-unimodal.
Label the vertices of a regular $n$-gon embedded in $\R^2$ in a clockwise fashion from $1$ to $n$.  
We define the \emph{circular distance} between two elements $i$ and $j$ of $\{1,2,\ldots,n\}$, denoted $\cd{i}{j}$, to be the number of edges in the shortest path between the vertices $i$ and $j$ of the $n$-gon.  
We also denote the path of shortest length from $i$ to $j$ by $\arc{i}{j}$.  
A subset $S\subset \{1,2,\ldots,n\}$ is called \emph{$r$-stable} if each pair $i,j\in S$ satisfies $\cd{i}{j}\geq r$.  
The \emph{$r$-stable $n,k$-hypersimplex}, denoted~$\Delta_{n,k}^{stab(r)}$, is the convex hull of the vectors $\sum_{i\in S}e_i$ where $S$ ranges over all $r$-stable $k$-subsets of the $n$-set.  
For fixed $n$ and $k$, these polytopes form the nested chain
\begin{equation*}\label{eqn:nested}
\Delta_{n,k}\supset\Delta_{n,k}^{stab(2)}\supset\Delta_{n,k}^{stab(3)}\supset\cdots\supset\Delta_{n,k}^{stab\left(\left\lfloor\frac{n}{k}\right\rfloor\right)} \, .
\end{equation*}
An approach to establishing unimodality for $h^*_{\Delta_{n,k}}$ is to attempt to inductively prove unimodality for $h^*_{\Delta_{n,k}^{stab(r)}}$ for all $r$, possibly using the following theorem due to Stanley.
\begin{theorem}[Stanley, \cite{stanleyhstarmonotonic}]
If $P\subseteq Q$ are lattice polytopes, then $h^*_{P,i}\leq h^*_{Q,i}$ for all~$i$, i.e. the $h^*$-vector for $Q$ dominates the $h^*$-vector for $P$ in each entry.
\end{theorem}
Jochemko and Sanyal~\cite{sanyaljochemko} have recently proved that this monotonicity property of $h^*$-vectors is equivalent to the non-negativity of $h^*$-vectors.

The following geometric properties of $\Delta_{n,k}^{stab(r)}$ are directly related to establishing $h^*$-unimodality, though they are also independently interesting.
\begin{theorem}[Braun and Solus, \cite{BraunSolusHyper}]
The regular unimodular triangulation of $\Delta_{n,k}$ given by the circuit triangulation restricts to a regular unimodular triangulation of $\Delta_{n,k}^{stab(r)}$ for all $r$.
Further, when $n$ is odd and $k=2$ there is a shelling of the circuit triangulation of $\Delta_{n,2}$ that first builds the $r$-stable hypersimplex and then builds the $(r-1)$-stable hypersimplex for every $1\leq r\leq \lfloor\frac{n}{2}\rfloor$.
\end{theorem}
\begin{theorem}[Hibi and Solus, \cite{hibisolushyper}]
Let $1<k<n-1$, and let $H$ denote the hyperplane in $\R^n$ defined by the equation $x_1+\cdots +x_n=k$.
For $1\leq r<\lfloor\frac{n}{k}\rfloor$, the facets of $\Delta_{n,k}^{stab(r)}$ are supported by the hyperplanes 
\[
\left\{x\in \R^n:x_\ell=0\right\}\cap H
\]
and 
\[
\left\{x\in \R^n:\sum_{i=\ell}^{\ell+r-1}x_i=1\right\}\cap H \, .
\]
Further, for $1\leq r< \lfloor\frac{n}{k}\rfloor$, the polytope $\Delta_{n,k}^{stab(r)}$ is Gorenstein if and only if $n=kr+k$.
\end{theorem}
The facet descriptions of the stable hypersimplices not provided here have also been obtained, see~\cite[Remark 2.3]{hibisolushyper}.

In \cite{BraunSolusHyper}, the $h^*$-polynomials of $\Delta_{n,2}^{stab(r)}$ are studied, and connections are obtained with Lucas and Fibonacci polynomials, independence polynomials of graphs, and monomial CR mappings of Lens spaces.
Since independence polynomials of graphs have been well-investigated with regard to unimodality, these connections yield additional tools which might be applied in this setting.
The vectors $h^*_{\Delta_{n,2}^{stab(r)}}$ are also determined to be unimodal for $r=2$, $r=3$, and $r=\lfloor\frac{n}{2}\rfloor -1$.
Many open questions about these polytopes remain, such as the following.
\begin{conjecture}[Braun and Solus, \cite{BraunSolusHyper}]
The $h^*$-vector of $\Delta_{n,2}^{stab(r)}$ is unimodal for all~$r$.
\end{conjecture}
\begin{question}[Braun and Solus, \cite{BraunSolusHyper}]
Is the $h^*$-vector of $\Delta_{n,k}^{stab(r)}$ unimodal for all~$r$?
\end{question}


\section{Log-Concavity and Real-Rootedness}

We conclude this survey with a brief discussion of log-concavity and real-rootedness for $h^*$-polynomials.
While a large amount of effort has focused on $h^*$-unimodality, the study of these two properties for lattice polytopes has also been a fruitful area of research.
We will discuss two such lines of investigation.


\subsection{$s$-Lecture Hall Polytopes}

Motivated by the study of lecture hall partitions~\cite{BME1}, given a finite sequence $s=(s_1,\ldots,s_n)$ Savage and Schuster~\cite{savageschuster} defined the $s$-lecture hall polytope to be
\[
P_n^s:=\left\{x\in\R^n:0\leq \frac{x_1}{s_1}\leq\frac{x_2}{s_2}\leq\cdots\leq\frac{x_n}{s_n}\leq 1\right\} \, .
\]
When $(s_1,s_2,\ldots,s_n)=(1,2,\ldots,n)$, this polytope contains a subset of the classical lecture hall partitions.
It is straightforward to verify that $P_n^s$ is a lattice polytope with vertex set
\[
\{(0,0,\ldots,0,s_k,s_{k+1},\ldots,s_n:1\leq k \leq n+1\} \, .
\]
Savage and Schuster derived a combinatorial interpretation of $h^*_{P_n^s}(z)$ using ascent statistics for $s$-inversion sequences.
Savage and Visontai~\cite{savagevisontai} used this combinatorial interpretation to prove the following.
\begin{theorem}[Savage and Visontai, \cite{savagevisontai}]
The polynomials $h^*_{P_n^s}(z)$ are real-rooted, hence the $h^*$-vector of $P_n^s$ is log-concave.
\end{theorem}
The proof given by Savage and Visontai uses only the combinatorial interpretation of these polynomials in terms of ascent sequences and \emph{compatible polynomial} techniques due to Chudnovsky and Seymour~\cite{chudnovskyseymour}.
\begin{question}
Is there a proof of real-rootedness, log-concavity, or unimodality for $h^*_{P_n^s}$ that fundamentally relies on the lattice point geometry of $\cone(P_n^s)$?
\end{question}
The following question regarding $P_n^s$ is natural in the context of unimodality.
\begin{question}
For which sequences $s$ is $P_n^s$ Gorenstein?
\end{question}

The lattice point geometry of $\cone(P_n^s)$ is related to the lattice point geometry of the $s$-lecture hall cone 
\[
C_n^s:=\left\{x\in \R^n:0\leq \frac{x_1}{s_1}\leq\frac{x_2}{s_2}\leq\cdots\leq\frac{x_n}{s_n}\right\} 
\]
for which the following is known.
\begin{theorem}[Beck, Braun, K\"{o}ppe, Savage, and Zafeirakopoulos, \cite{sLectHallGorenstein}]
Let $\ell > 0$ and $b \not = 0$ be integers satisfying $\ell^2+4b \geq 0$.
Let $s=(s_1, s_2, \ldots)$ be defined by
\[
s_n = \ell s_{n-1} + b s_{n-2} \, ,
\]
with initial conditions $s_1=1$, $s_0=0$.
For each $n$, let $C_n^s$ be defined by the truncation $(s_1,\ldots,s_n)$.
Then  $C_n^s$ is Gorenstein for all $n \geq 0$ if and only if $b=-1$.
If $b \not = -1$, there exists $n_0 = n_0(b,\ell)$ such that 
$C_n^s$
 fails to be Gorenstein
for all $n \geq n_0$.
\end{theorem}
For example, when $\ell=2$ and $b=-1$, we obtain the sequence $(1,2,3,4,\ldots)$, corresponding to the classical lecture hall partition cone.


\subsection{Dilations of Polytopes and Multivariate Techniques}

Our final topic regards the behavior of $h^*$-vectors for integral dilates $nP$ of a lattice polytope $P$.
This work has generally been motivated by the following theorem.
\begin{theorem}[Knudsen, Mumford, and Waterman, \cite{knudsenmumford}]
For every lattice polytope $P$, there exists an integer $n$ such that $nP$ admits a regular unimodular triangulation.
\end{theorem}
Identifying the values of $n$ that correspond to $nP$ with regular unimodular triangulations is a deep and subtle open problem.
Bruns and Gubeladze~\cite{brunsgubeladzeunicover} proved that for each $d$, there exists a natural number $c_d$ such that for all $c>c_d$ and all $d$-dimensional lattice polytopes $P$, $cP$ is equal to the union of the unimodular subsimplices of $cP$.
It is unknown whether this is true for the property of $cP$ admitting a regular unimodular triangulation.

Recently there has been interest in the study of the transformations taking a power series
\[
\sum_{n\geq 0}a_nz^n = \frac{h(z)}{(1-z)^d}
\]
to
\[
\sum_{n\geq 0}a_{rn}z^n = \frac{h^{<r>}(z)}{(1-z)^d}
\]
for $r\geq 1$.
For a lattice polytope $P$, this corresponds to studying $h^*$-polynomials for integer dilates of $P$. 
The following result strengthens a previous theorem due to Brenti and Welker~\cite{brentiwelkerveronese} regarding the Veronese construction for graded algebras and associated Hecke operators.
\begin{theorem}[Beck and Stapledon, \cite{beckstapledon}]
Fix a positive integer $d$ and let $\rho_1<\rho_2<\cdots <\rho_d=0$ denote the roots of $zA_d(z)$, where $A_d(z)$ is the Eulerian polynomial.
There exists positive integers $m_d$ and $n_d$ such that, if $P$ is a $d$-dimensional lattice polytope and $n\geq n_d$, then $h^*_{nP}(z)$ has negative real roots $\beta_1(n)<\beta_2(n)<\cdots <\beta_{d-1}(n)<\beta_d(n)<0$ with $\beta_i(n)\to\rho_i$ as $n\to \infty$, and the coefficients $\left\{h^*_i(n)\right\}_{1\leq i\leq d}$ of $h^*_{nP}(z)$ are positive, strictly log-concave, and satisfy
\begin{align*}
1=h^*_0(n)<h^*_d(n)<h^*_1(n)<\cdots \hspace{4cm}\\
<h^*_i(n)<h^*_{d-i}(n) <h^*_{i+1}(n)<\cdots \hspace{2cm} \\
<h^*_{\left\lfloor \frac{d+1}{2}\right\rfloor}(n)<m_dh^*_d(n) \, .
\end{align*}
\end{theorem}
A strong conjecture regarding real-rootedness for integral dilates of lattice polytopes is the following.
\begin{conjecture}[Beck and Stapledon, \cite{beckstapledon}]
If $P$ is $d$-dimensional, then for $n\geq d$ the polynomial $h^*_{nP}(z)$ has distinct, negative real roots.
\end{conjecture}
In support of this conjecture, Higashitani~\cite{higashitanipersonal} has recently proved that for a lattice polytope $P$, if $n\geq \deg(h^*_P(z))$ then the $h^*$-vector of $nP$ is strictly log-concave.

McCabe and Smith~\cite{mccabesmith} extended these investigations to the setting of multivariate log-concavity for multigraded Hilbert series, and introduce new techniques for establishing log-concavity.
They use these techniques to provide a new proof of the log-concavity of the Eulerian polynomials.
The study of log-concavity and unimodality for multivariate $h^*$-polynomials deserves further investigation, see~\cite[Section 3]{mccabesmith} for further discussion and open problems.

\begin{remark}
Another multivariate approach to Ehrhart $h^*$-polynomials are \emph{local} $h^*$-polynomials.
Katz and Stapledon~\cite[Theorem 9.4]{katzstapledonlocal} have recently established various unimodality results for these objects.
\end{remark}


\begin{acknowledgement}

The author is partially supported by the National Security Agency through award H98230-13-1-0240.
Thanks to Christos Athanasiadis, Matthias Beck, Robert Davis, Jesus De Loera, David Haws, Takayuki Hibi, Akihiro Higashitani, Carla Savage, and Liam Solus for their helpful comments and suggestions.

\end{acknowledgement}


\bibliographystyle{spmpsci}
\bibliography{Braun}

\newpage 

\section*{Appendix}

Consider the following three properties of a polynomial $p=a_0+a_1z+\cdots+a_dz^d$ with $a_d\neq 0$ and each $a_i\geq 0$:
\begin{enumerate}
\item $p$ is called \emph{unimodal} if there exists an index $j$ such that $a_i\leq a_{i+1}$ for all $i<j$ and $a_{i}\geq a_{i+1}$ for all $i\geq j$.  
\item $p$ is called \emph{log-concave} if $a_j^2\geq a_{j-1}a_{j+1}$ for all $1\leq j\leq d-1$.
\item $p$ is called \emph{flawless and non-decreasing in the first half} if 
\begin{itemize}
\item $a_i\leq a_{d-i}$ for every $0\leq i\leq \lfloor d/2\rfloor$, and
\item $a_0\leq a_1\leq \cdots \leq a_{\lfloor d/2\rfloor}$.
\end{itemize}
\end{enumerate}

An open problem in commutative algebra is whether or not the $h$-polynomial of a standard graded Cohen-Macaulay integral domain is unimodal.
In the published version of this survey, it is incorrectly stated after Conjecture~1 that counterexamples exist to unimodality in this setting.
The examples referenced are those found in~\cite{MR1230867} and~\cite{MR1361577}, which provide counterexamples to the conditions of being (A) log-concave and (B) flawless and non-decreasing in the first half.

Thus, the correct statement is: not every standard graded Cohen-Macaulay integral domain has an $h$-polynomial that is (A) log-concave or (B) flawless and non-decreasing in the first half.
To the knowledge of the author it is currently unknown whether or not such rings always have unimodal $h$-polynomials.

\end{document}